\documentclass[11pt]{amsart}
\usepackage{amsfonts}
\usepackage{mathrsfs}
\usepackage{amssymb,latexsym}
\usepackage{pstcol,pstricks,color}

 \numberwithin{equation}{section}

\title{Some Identities Involving Three Kinds of Counting Numbers}

\begin{document}
\maketitle

\begin{center}
L. C. Hsu

Mathematics Institute, Dalian University of Technology, Dalian 116024, China\\[5pt]

\end{center}\vskip0.5cm

\subsection*{Abstract} In this note, we present several identities involving
binomial coefficients and the two kind of Stirling numbers.

\vskip0.5cm

{\large\bf 1. Introduction}

\vspace{12pt}

Adopting Knuth's notation, let us denote by
$\left[\begin{array}{c}n\\k\end{array}\right]$ and
$\left\{\begin{array}{c}n\\k\end{array}\right\}$ the unsigned
(absolute) Stirling number of the first kind and the ordinary
Stirling number of the second kind, respectively. Here in
particular,
$\left[\begin{array}{c}0\\0\end{array}\right]=\left\{\begin{array}{c}0\\0\end{array}\right\}=1$,
$\left[\begin{array}{c}n\\0\end{array}\right]=\left\{\begin{array}{c}n\\0\end{array}\right\}=0~(n>0)$,
and
$\left[\begin{array}{c}n\\k\end{array}\right]=\left\{\begin{array}{c}n\\k\end{array}\right\}=0~(0\le
n<k)$. Generally, $\left[\begin{array}{c}n\\k\end{array}\right]$,
$\left\{\begin{array}{c}n\\k\end{array}\right\}$ and the binomial
coefficients $\left(\begin{array}{c}n\\k\end{array}\right)$ may be
regarded as the most important counting numbers in combinatorics.
The object of this short note is to propose some combinatorial
identities each consisting of these three kinds of counting numbers,
namely the following
$$
\sum_k\left[\begin{array}{c}k\\p\end{array}\right]\left\{\begin{array}{c}n+1\\k+1\end{array}\right\}(-1)^k=
\left(\begin{array}{c}n\\p\end{array}\right)(-1)^p
$$
$$
\sum_k\left[\begin{array}{c}k+1\\p+1\end{array}\right]\left\{\begin{array}{c}n\\k\end{array}\right\}(-1)^k=
\left(\begin{array}{c}n\\p\end{array}\right)(-1)^n
$$
$$
\sum_{j,k}\left[\begin{array}{c}n\\k\end{array}\right]\left\{\begin{array}{c}k\\j\end{array}\right\}
\left(\begin{array}{c}n\\j\end{array}\right)(-1)^k=(-1)^n
$$
$$
\sum_{j,k}\left\{\begin{array}{c}n\\k\end{array}\right\}\left[\begin{array}{c}k\\j\end{array}\right]
\left(\begin{array}{c}n\\j\end{array}\right)(-1)^k=(-1)^n
$$
$$
\sum_{j,k}\left(\begin{array}{c}n\\k\end{array}\right)\left\{\begin{array}{c}k\\j\end{array}\right\}
\left[\begin{array}{c}j+1\\p\end{array}\right](-1)^j=\left\{\begin{array}{cl}0,
&(n+1>p)\\[6pt]
(-1)^n, &(n+1=p)\end{array}\right.
$$
$$
\sum_{j,k}\left[\begin{array}{c}n\\k\end{array}\right]\left(\begin{array}{c}k\\j\end{array}\right)
\left\{\begin{array}{c}j+1\\p\end{array}\right\}(-1)^j=\left\{\begin{array}{cl}0,
&(n+1>p)\\[6pt]
(-1)^n, &(n+1=p)\end{array}\right.
$$
Here each of the summations in (1) and (2) extends over all $k$ such
that $0\le k\le n$ or $p\le k\le n$, and all the double summations
within (3)--(6) are taken over all possible integers $j$ and $k$
such that $0\le j\le k\le n$.

Note that (1) is a well-known identity that has appeared in the
Table 6.4 of Graham-Knuth-Patashnik's book [1] (cf. formula (6.24)).
It is quite believable that (1) and (2) may be the most simple
identities each connecting with the three kinds of counting numbers.
\\[10pt]
{\large\bf 2. Proof of the identities}

\vspace{10pt}

In order to verify (2)--(6), let us recall that the orthogonality
relations
$$
\sum_{k}\left[\begin{array}{c}n\\k\end{array}\right]\left\{\begin{array}{c}k\\p\end{array}\right\}
(-1)^{n-k}=\sum_k\left\{\begin{array}{c}n\\k\end{array}\right\}\left[\begin{array}{c}k\\p\end{array}\right]
(-1)^{n-k}=\delta_{np}=\left\{\begin{array}{ll}0, &(n\neq p)\\[8pt]
1, &(n=p)\end{array}\right.
$$
are equivalent to the inverse relations
$$
a_n=\sum_{k}\left[\begin{array}{c}n\\k\end{array}\right](-1)^{n-k}b_k
\Leftrightarrow
b_n=\sum_k\left\{\begin{array}{c}n\\k\end{array}\right\}a_k.
$$
Also, we shall make use of two known identities displayed in the
Table 6.4 of [1], viz.
$$
\left\{\begin{array}{c}n+1\\p+1\end{array}\right\}=\sum_k\left(\begin{array}{c}n\\k\end{array}\right)
\left\{\begin{array}{c}k\\p\end{array}\right\}
$$
$$
\left[\begin{array}{c}n+1\\p+1\end{array}\right]=\sum_k\left[\begin{array}{c}n\\k\end{array}\right]
\left(\begin{array}{c}k\\p\end{array}\right).
$$
Now, take
$a_n=(-1)^n\left[\begin{array}{c}n+1\\p+1\end{array}\right]$ and
$b_k=(-1)^k\left(\begin{array}{c}k\\p\end{array}\right)$, so that
(10) can be embedded in the first equation of (8). Thus it is seen
that (10) can be inverted via (8) to get the identity (2).

(3) and (4) are trivial consequences of (7). Indeed rewriting (7) in
the form

\vspace{8pt} ~~~~~~~~~~~~~~~~~$\displaystyle
\sum_k\left[\begin{array}{c}n\\k\end{array}\right]\left\{\begin{array}{c}k\\j\end{array}\right\}(-1)^k=
\sum_k\left\{\begin{array}{c}n\\k\end{array}\right\}\left[\begin{array}{c}k\\j\end{array}\right](-1)^k=
(-1)^n\delta_{nj}$\hfill $(7)'$
\\[8pt]
and noticing that
$\sum_j\left(\begin{array}{c}n\\j\end{array}\right)\delta_{nj}=
\left(\begin{array}{c}n\\n\end{array}\right)=1$, we see that
(3)--(4) follow at once from $(7)'$.

For proving (5), let us make use of (9) and (7) with $p$ being
replaced by $j$. We find
$$
\sum_j\left\{\begin{array}{c}n+1\\j+1\end{array}\right\}\left[\begin{array}{c}j+1\\p\end{array}\right](-1)^j=
\sum_j\sum_k\left(\begin{array}{c}n\\k\end{array}\right)\left\{\begin{array}{c}k\\j\end{array}\right\}
\left[\begin{array}{c}j+1\\p\end{array}\right](-1)^j=(-1)^n\delta_{n+1,p}.
$$
Hence (5) is obtained.

Similarly, (6) is easily derived from (10) and (7).
\\[10pt]
{\large\bf 3. Questions}

\vspace{10pt}

It may be a question of certain interest to ask whether some of the
identities (1)--(6) could be given some combinatorial
interpretations with the aid of the inclusion-exclusion principle or
the method of bijections. Also, we have not yet decided whether
(1)--(6) could be proved by the method of generating functions (cf.
[2]).


{\bf AMS Classification Numbers}: 05A10, 05A15, 05A19.

\end{document}